\documentclass[12pt,twoside]{article}
\usepackage[centertags]{amsmath}
\usepackage{amsfonts}
\usepackage{amssymb}
\usepackage{amsthm}
\usepackage{newlfont}
\usepackage{color}

\date{\color{green}  2013 September 16}
\def \n {\noindent}
\usepackage{fancyhdr}
\begin{document}

\centerline{{\bf \color{red} 2013 September 16}}
\centerline{}
\centerline{}
\centerline {\Large{\bf  On chaoticity of the sum of chaotic shifts
with their adjoints }}
\centerline{}
\centerline {\Large{\bf in Hilbert space and applications}}
\centerline{}
\centerline {\Large{\bf to some chaotic weighted shifts acting}}
\centerline{}
\centerline {\Large{\bf on some Fock-Bargmann spaces}}
\centerline{}
\centerline{ Abdelkader  Intissar}
\centerline{}
\centerline{\it Equipe d'Analyse spectrale,Facult\'e des Sciences et Techniques}
\centerline{ Universit\'e de Cort\'e, 20250 Cort\'e, France}
\centerline{\it T\'el: 00 33 (0) 4 95 45 00 33-\it Fax: 00 33 (0) 4 95 45 00 33}
\centerline{\it intissar@univ-corse.fr}
\centerline{\&}
\centerline{\it Le Prador, 129, rue du Commandant Rolland, 13008 Marseille,
France}
\centerline{}
\newtheorem{Theorem}{\quad Theorem}[section]
\newtheorem{Definition}[Theorem]{\quad Definition}
\newtheorem{Corollary}[Theorem]{\quad Corollary}
\newtheorem{Lemma}[Theorem]{\quad Lemma}
\newtheorem{Example}[Theorem]{\quad Example}
\newtheorem{Remark}[Theorem]{\quad Remark}
\begin{abstract} \small{This article is intended to outline some the recent work by the author on the chaoticity of some specific bakward shift unbounded operators realized as differential operators acting on some Fock- Bargmann spaces and give sufficient conditions on a linear unbounded densely defined chaotic shift operator $\mathbb{T}$ acting on a Hilbert space for the operator $\mathbb{T} + \mathbb{T}^{*}$ to be chaotic where $\mathbb{T}^{*}$ is its adjoint.\\}.
\end{abstract}
\newpage
\section{ Introduction}
A continuous operator $\mathbb{T}$ on a Banach space $\mathbb{X}$ is said to be hypercyclic if the following condition is met:\\
There exists an element $\phi \in \mathbb{X}$ that its orbit $Orb(\mathbb{T},\phi) = \{\phi, \mathbb{T}\phi, \mathbb{T}^{2}\phi, .....\}$ is dense in $\mathbb{X}$ and is said to be chaotic in the sense of Devaney {\bf [2,15]} if the following conditions is met:\\
1) $\mathbb{T}$ is hypercyclic.\\
2) The set $\{\phi \in \mathbb{X}; \exists \quad n \in I\!\!N$ such that $\mathbb{T}^{n}\phi = \phi\}$ of periodic points of operator $\mathbb{T}$ is dense in $\mathbb{X}$.\\
It is well known that linear operators in finite-dimensional linear spaces can't be chaotic but the nonlinear operator may be. Only in infinite-dimensional linear spaces can linear operators have chaotic properties.These last properties are based on the phenomenon of hypercyclicity or the phenomen of non-wandercity.\\
The study of the phenomenon of hypercyclicity originates in the papers by Birkoff {\bf[7]} and Maclane {\bf [27]} that show, respectively, that the operators of translation and differentiation, acting on the space of entire functions are hypercyclic. \\
The theories of hypercyclic operators and chaotic operators have been intensively developed for bounded linear operator, we refer to {\bf [13,14]} and references therein and for a bounded operator, Ansari  asserts in  {\bf[1]} that powers of a hypercyclic bounded operator are also hypercyclic\\

 For an unbounded operator, Salas exhibit in {\bf [30]} an unbounded hypercyclic operator whose square is not hypercyclic. The result of Salas show that one must be careful in the formal manipulation of operators with restricted domains. For such operators it is often more convenient to work with vectors rather than with operators themselves.\\
 Now, let $\mathbb{T}$ be an unbounded operator on a separable infinite dimensional Banach space $\mathbb{X}$.\\

\noindent We define the following sets:\\

 \noindent $D(\mathbb{T}) = \{\phi \in \mathbb{X} ; \mathbb{T}\phi \in \mathbb{X}\}$ $\hfill { }(1.1)$\\

 \noindent $D(\mathbb{T}^{\infty}) = \bigcap_{n=0}^{\infty}D(\mathbb{T}^{n})$ $\hfill { }(1.2)$\\

 The notion of chaos for unbounded operators was defined in {\bf [6]} by B\'es et al as follows:\\

\begin{Definition}
 {\it A linear unbounded densely defined operator $(\mathbb{T}, D(\mathbb{T}))$ on a Banach space $\mathbb{X}$ is called chaotic if the following conditions are met:\\
 1) $\mathbb{T}^{n }$ is closed for all positive integers $n$..\\
 2) there exists an element $\phi \in D(\mathbb{T}^{\infty})$ whose orbit $Orb(\mathbb{T},\phi) = \{\phi, \mathbb{T}\phi, \mathbb{T}^{2}\phi, .....\}$ is dense in $\mathbb{X}$\\
 3) the set $\{\phi \in \mathbb{X}; \exists \quad m \in I\!\!N$ such that $\mathbb{T}^{m}\phi = \phi\}$ of periodic points of operator $\mathbb{T}$ is dense in $\mathbb{X}$.}\\
 \end{Definition}

Recently these theories are begin developed on some concrete examples of unbounded linear operators, see {\bf [5,8,16]}. In {\bf [16]} it has been shown that the operators $\mathbb{H}_{p} = z^{p}\frac{d^{p+1}}{dz^{p+1}} ; p=0, 1, .....$ are chaotic in the sense of Definition 1.1 on the classic Bargmann space {\bf [3]} of entire functions with $e^{-\mid z\mid^{2}}$ measure.\\

In {\bf [17]} we have considered generalized Bargmann spaces (the spaces of entire functions with $e^{-\mid z\mid^{\beta}}$ measure; $\beta > 0$) and we have proved that the operators $\mathbb{H}_{p} = z^{p}\mathbb{D}^{p+1} ; p=0, 1, .....$ in these spaces are chaotic where  $\mathbb{D}$ is the adjoint operator of the operator of multiplication by the independent variable $z$  on these spaces. $\mathbb{D}$ belongs to class Gelfond-Leontiev operators of generalized differentiation {\bf [10]}\\

In {\bf [18]} we have considered non-compact lattice of $(\Gamma,\chi)$-theta Fock-Bargmann spaces  and we have proved that the operators $\mathbb{H}_{p} = e^{2ipz}\mathbb{D}^{p+1} ; p=0, 1, .....$ in these spaces are chaotic where  $\mathbb{D}$ is the adjoint operator of the operator of multiplication by the function $M(z) = e^{2ipz}$  on these spaces.\\

In the present work, we give sufficient conditions on a linear unbounded densely defined chaotic shift operator $\mathbb{T}$ acting on a Hilbert space such that $\mathbb{T} + \mathbb{T}^{*}$ is chaotic where $\mathbb{T}^{*}$ is its adjoint and we apply this sufficient conditions to above operators.\\
In appendix, we consider Fock-Bargmann space on the unit disk quantized where we will consider the annihilation operator $\mathbb{A}$ associated its orthonormal basis and we prove that the operators $\mathbb{H}_{p} = \mathbb{A}^{*^{p}}\mathbb{A}^{p+1} ; p=0, 1, .....$ are chaotic where  $\mathbb{A}^{*}$ is the adjoint operator of the annihilation operator $\mathbb{A}$.\\

This paper is organized as follows :\\
In section 2 we recall some sufficient conditions on hypercyclicity of bounded operators given by Godefroy-Shapiro's lemma {\bf [12]} or on hyperccylicity of unbounded operators given by B$\grave{e}$s-Chan-Seubert theorem {\bf [6]}.\\
We recall also some elementary properties of classic Bargamann space, the generalized Bargmann space, the $(\Gamma,\chi)$-theta Fock-Bargmann space, Fock-Bargmann space on the unit disk quantized and the action of $\mathbb{H}_{p}; p \in \mathbb{N}$ on these spaces.\\
In section 3 we give sufficient conditions on a linear unbounded densely defined chaotic shift operator $\mathbb{T}$ acting on a Hilbert space such that $\mathbb{T} + \mathbb{T}^{*}$ is chaotic where $\mathbb{T}^{*}$ is its adjoint with application to $\mathbb{H}_{p} ; p \in \mathbb{N}$ defined on above spaces.\\
In appendix we consider the annihilation operator $\displaystyle{\mathbb{A} = \frac{d}{dz}}$ acting on orthonormal basis of the Fock-Bargmann space on the unit disk quantized. We prove that the operators $\mathbb{H}_{p} = \mathbb{A}^{*^{p}}\mathbb{A}^{p+1} ; p=0, 1, .....$ acting on this space are chaotic where  $\displaystyle{\mathbb{A}^{*}= z^{2}\frac{d}{dz} + 2\nu z; \nu > 1}$ is the adjoint operator of $\mathbb{A}$.\\
As these operators are unilateral weighted backward shifts with an explicit weight, we use the results of B$\grave{e}$s et al to proof the chaoticity of $\mathbb{H}_{p}$ on Fock-Bargmann space associated to Poincar\'e disk (we can also use the results of Bermudez et al {\bf [5]} to proof the chaoticity of our operators $\mathbb{H}_{p}$).\\

\section{Action of operators $\mathbb{H}_{p}$ of order $p$ on associated Fock-Bargmann spaces}

Before to recall the Fock-Bargmann spaces those the operators $\mathbb{H}_{p}$ with  domain $D(\mathbb{H}_{p} )$ acting, we begin by to recall that an unbounded operator $\mathbb{T}$ is hypercyclic if there is a vector $\phi$ in the domain of $\mathbb{T}$ such that for every integer $m > 1$ the vector $T^{m}\phi$ is in the domain of $\mathbb{T}$ and the orbit $\{\phi, \mathbb{T}\phi, \mathbb{T}^{2}\phi, \mathbb{T}^{3}\phi,... \}$ is dense in $\mathbb{X}$.\\

We recall also the Godefroy-Shapiro's lemma {\bf[12]} and B$\grave{e}$s-Chan-Seubert's theorem {\bf[6]} that we will used in section 3 and in appendix\\

A) hypercyclicity criterion of Godefroy-Shapiro and of B$\grave{e}$s-Chan-Seubert\\

\begin{Lemma} (Godefroy-Shapiro ({[12]}).\\
 {\it Let $\mathbb{X}$ be a separable Fr\'echet space and $\mathbb{T}$ is bounded operator on $\mathbb{X}$ and $\mathbb{Y}_{1}$, $\mathbb{Y}_{2}$ are two dense subsets of $\mathbb{X}$ and $\mathbb{S} : \mathbb{Y}_{1} \rightarrow \mathbb{Y}_{1}$ such that :\\
(1) $\mathbb{T}\mathbb{S}\phi = \phi, \quad \forall \quad \phi \in \mathbb{Y}_{1}$ \\
(2) $lim \mathbb{S}^{m}\phi = 0, \quad \forall \quad \phi \in Y_{1}$ as $m \rightarrow +\infty$\\
(3) $lim \mathbb{T}^{m}\phi = 0, \quad \forall \quad \phi \in Y_{2}$ as $m \rightarrow +\infty$\\
then $\mathbb{T}$ is hypercyclic operator.}\\
\end{Lemma}

then hypercyclic vectors can be constructed for an unbounded operator $\mathbb{T}$, under a sufficient condition analogous to the above hypercyclicity criterion witch is given by $B\grave{e}s-Chan-Seubert$'s theorem\\

 \begin{Theorem}
 (B$\grave{e}$s-Chan-Seubert  {\bf [6]}, p.258 )\\{\it Let $\mathbb{X}$ be  a separable infinite dimensional Banach and let $\mathbb{T}$ be a densely defined linear operator on $\mathbb{X}$. Then $\mathbb{T}$ is hypercyclic if\\
(i) $\mathbb{T}^{m}$ is closed operator for all positive integers $m$.\\
(ii) There exist a dense subset $\mathbb{Y}$ of the domain $D(\mathbb{T})$ of $\mathbb{T}$ and a (possibly nonlinear and discontinuous) mapping $\mathbb{S} : \mathbb{Y} \longrightarrow \mathbb{Y} $ so that $\mathbb{T}\mathbb{S} = I_{\mid \mathbb{Y}}$ ($I_{\mid \mathbb{Y}}$ is identity on $\mathbb{Y}$) and $ \mathbb{T}^{n}, \mathbb{S}^{n} \longrightarrow 0 $ pointwise on $\mathbb{Y}$ as $ n \longrightarrow \infty.$}\\
\end{Theorem}

\noindent B) the classic Bargmann space {\bf [3]} is defined by:\\

  $\mathbb{B} = \{\phi: I\!\!\!\!C\longrightarrow  I\!\!\!\!C\, entire ; \displaystyle {\int_{I\!\!\!\!C}}\displaystyle{\mid\phi(z)\mid^{2}}e^{-\mid z\mid^{2}}dxdy < \infty \}$ $\hfill { }(2.1)$\\
\quad where $z = x+iy$.\\

$\mathbb{B}$ is a Hilbert space with an inner product\\

$< \phi,\psi > = \displaystyle {\int_{I\!\!\!\!C}}\displaystyle{\phi(z)\overline{\psi(z)}e^{-\mid z\mid^{2}}dxdy}$$\hfill { }(2.2)$\\

\noindent and the associated norm is denoted by $\mid\mid . \mid\mid $.\\

The functions $e_{n}(z) = \frac{z^{n}}{\sqrt{n!}}; n = 0, 1, 2, ....$ form a complete orthonormal set in $\mathbb{B}$\\

The operator of multiplication by the independent variable $z$  on $\mathbb{B}$ is defined by :\\

\noindent $\mathbb{A}^{*}\phi(z) = z\phi(z)$ with domain $D(\mathbb{A}^{*}) = \{\phi \in \mathbb{B}; z\phi \in \mathbb{B}\}\hfill { } (2.3)$\\

The operator $\mathbb{A}^{*}$ acts on $e_{n}(z)$ as following:\\

$\mathbb{A}^{*}e_{n}(z) = \sqrt{n+1}e_{n+1}(z) = \omega_{n}e_{n+1}(z)$ with $\omega_{n} = \sqrt{n+1}$$\hfill { }(2.4)$\\

Then its adjoint is differentiation operator $\displaystyle{\mathbb{A}\phi(z) = \frac{d}{dz}\phi(z)}$ with domain $D(\mathbb{A}) = \{\phi \in \mathbb{B}; \frac{d}{dz}\phi \in \mathbb{B}\}\hfill { }(2.5) $\\

it is given also by :\\

$\mathbb{A}e_{0}(z) = 0$ and
$\mathbb{A}e_{n}(z) = \sqrt{n}e_{n-1}(z) =  \omega_{n-1}e_{n-1}(z), n \geq 1 $ $\hfill { }(2.6)$\\

We define now a family of weighted shifts $\mathbb{H}_{p}$ acting on $\mathbb{B}$ as following\\

$\mathbb{H}_{p} = \mathbb{A}^{*^{p}}\mathbb{A}^{p+1}$ with domain $D(\mathbb{H}_{p}) = \{\phi \in \mathbb{B}; \mathbb{H}_{p}\phi \in \mathbb{B}\}$ $\hfill { }(2.7)$\\

\noindent Then we get\\

$\mathbb{H}_{p}^{*}e_{n}(z) = \mathbb{A}^{*^{p+1}}\mathbb{A}^{p}e_{n}(z) =$\\

$\sqrt{(n+1)}\prod_{j=1}^{p}(n-j)e_{n+1}(z)$ for $ n \ge p \geq 0 $ $\hfill { }(2.8)$\\

\noindent i.e. $ \mathbb{H}_{p} ^{*}$ is weighted shift with weight \\

$ \omega_{n,p} = \omega_{n}\prod_{j=1}^{p}\omega_{n-j}^{2}$ for $ n \ge p \geq 0 $ $\hfill { }(2.9)$\\\\

\noindent C) We define the generalized Bargmann space by :\\

$\mathfrak{F}_{\beta} = \{\phi: I\!\!\!\!C\longrightarrow  I\!\!\!\!C\, entire ; \displaystyle {\int_{I\!\!\!\!C}}\displaystyle{\mid\phi(z)\mid^{2}}e^{-\mid z\mid^{\beta}}d\mu(z) < \infty \}$ $\hfill { }(2.10)$\\
where $\beta > 0$ is an arbitrary constant, $d\mu(z) = \frac{\beta}{2\pi \Gamma(\frac{2}{\beta})} dxdy $ and $z = x+iy$.\\

Note that $\mathfrak{F}_{2}$ coincides with the classic Bargmann space.\\

$\mathfrak{F}_{\beta}$ is a Hilbert space with an inner product\\

$< \phi,\psi > = \frac{\beta}{2\pi\Gamma(\frac{2}{\beta})}\displaystyle {\int_{I\!\!\!\!C}}\displaystyle{\phi(z)\overline{\psi(z)}e^{-\mid z\mid^{\beta}}dxdy}\hfill { }(2.11)$\\

\noindent and the associated norm is denoted by $\mid\mid . \mid\mid $.\\

Let $m_{0} = 0 $, $m_{n} = \frac{\Gamma(\frac{2}{\beta}(n+1))}{\Gamma(\frac{2n}{\beta})}$
$n =1,2, ...$ and $[m_{n}]! = m_{1}.m_{2}......m_{n}$ then it may be shown that the functions\\

$e_{0}(z) = 1$ and $e_{n}(z) = \frac{z^{n}}{\sqrt{[m_{n}]!}}; n = 1, 2, ....\hfill { }(2.12)$\\

form a complete orthonormal set in $\mathfrak{F}_{\beta}$.\\

Define the principal vectors $e_{\lambda} \in \mathfrak{F}_{\beta}$ (for every $\lambda \in I\!\!\!\!C $) as complex valued functions \\

$ e_{\lambda}(z) = e(z,\lambda) $= $1 + \displaystyle{\sum_{n=1}^{\infty}e_{n}(z)\overline{e_{n}(\lambda)}}$ of  $\lambda$ and $z$ in $I\!\!\!\!C$ \\

\noindent If $\phi(z) = \displaystyle{\sum_{n}^{\infty}a_{n}e_{n}(z)}$ then
$ < \phi,  e_{\lambda} > = \phi(\lambda)$ \quad \quad ( the reproducing property)
\noindent because
$\displaystyle {\int_{I\!\!\!\!C} \sum_{n}^{\infty}a_{n}e_{n}(z)(\overline{1 + \sum_{n=1}^{\infty}e_{n}(z)e_{n}(\lambda)}) e^{-\mid z\mid^{\beta}}d\mu(z)}$ = $\displaystyle {a_{0} +\sum_{n=1}^{\infty}a_{n}e_{n}(\lambda)\mid\mid e_{n}\mid\mid = \phi(\lambda)}$\\
\noindent or, in other words\\

$\phi (z) = \displaystyle {\int_{I\!\!\!\!C}}\displaystyle{\phi(\lambda)\overline{e_{\lambda}(z)}}e^{-\mid \lambda\mid^{\beta}}d\mu(\lambda)$ for all $\phi \in \mathfrak{F}_{\beta}$ $\hfill { }(2.13)$\\

\noindent so that  $e_{\lambda}(z)$ is called a reproducing kernel for $\mathfrak{F}_{\beta}$\\

Note that the reproducing kernel $ e_{\lambda}(z)$ is uniquely determined by the
Hilbert space $\mathfrak{F}_{\beta}$ and the evaluation linear functional $\phi \in \mathfrak{F}_{\beta} \rightarrow \phi(z) \in I\!\!\!\!C$ is a bounded linear functional on $\mathfrak{F}_{\beta}$.\\

So applying (2.13) to the function $e_{z}$ at $\lambda$; we get $e_{z}(\lambda) = < e_{z};e_{\lambda} >$ \quad \quad for $z$; $\lambda \in I\!\!\!\!C$ and by the above relations, for $ z \in I\!\!\!\!C$ we obtain \\

$\mid\mid e_{z}\mid\mid = \sqrt{< e_{z}, e_{z} >} = \sqrt{e(z, z)}$. $\hfill { }(2.14)$\\

A systematic study of these generalized Bargmann spaces can be founded in {\bf [24] } where Irac-Astaud and Rideau have constructed an deformed harmonic algebra (DHOA) on $\mathfrak{F}_{\beta}$ and in {\bf [25] } where Knirsch and Schneider have invesigated the continuity and Schatten–von Neumann $p$-class membership of Hankel operators with anti-holomorphic symbols on these spaces with $\beta \in \mathbb{N}$.\\
Note that  the generalized Bargmann spaces $\mathfrak{F}_{\beta}$ are different from the generalized Bargmann spaces $\mathbb{E}_{m}$ $ m= 0, 1, ....$ defined in {\bf [19]}. It would be interesting to characterize the orthogonal space of $\mathfrak{F}_{\beta}$ in $\mathbb{L}_{2}(\mathbb{C}, e^{-\mid z\mid^{\beta}}d\mu(z))$ for $\beta \neq 2$.\\

On the generalized Bargmann representation $\mathfrak{F}_{\beta}$, we denote now the operator of multiplication by the independent variable $z$  on $ \mathfrak{F}_{\beta}$ by :\\

\noindent $\mathbb{M}\phi(z) = z\phi(z)$ with domain $I\!\!D(\mathbb{M}) = \{\phi \in \mathfrak{F}_{\beta}; z\phi \in \mathfrak{F}_{\beta}\}\hfill { } (2.15)$\\

The operator $\mathbb{M}$ acts on $e_{n}(z)$ as following:\\

$\mathbb{M}e_{n}(z) = \frac{\sqrt{\Gamma(\frac{2}{\beta}(n+2))}}{\sqrt{\Gamma(\frac{2}{\beta}(n+1))}}e_{n+1}(z) \hfill { } (2.16)$\\

Then its adjoint is generalized differentiation given by :\\

$\mathbb{D}e_{n}(z) = \frac{\sqrt{\Gamma(\frac{2}{\beta}(n+1))}}{\sqrt{\Gamma(\frac{2}{\beta}n)}}e_{n-1}(z) \hfill { } (2.17)$\\

\noindent and for $\phi(z) = \displaystyle{\sum_{n=0}^{\infty}a_{n}z^{n}}$ we have
$\mathbb{D}1 = 0$ and $\mathbb{D}\phi(z) = \displaystyle{\frac{1}{z}\sum_{n=0}^{\infty}a_{n}m_{n}z^{n}}$
where $ m_{n} = \frac{\Gamma(\frac{2}{\beta}(n+1))}{\Gamma(\frac{2}{\beta}n)}$
with domain:\\

$D(\mathbb{D}) = \{\phi \in \mathfrak{F}_{\beta}; \mathbb{D}\phi \in \mathfrak{F}_{\beta}\} \hfill { } (2.18)$\\

Note that if $\beta = 2$ the generalized differentiation operator $\mathbb{D}$ is:\\

$\displaystyle{\mathbb{D}\phi(z) = \frac{d}{dz}\phi(z)}$  $\hfill { } (2.19)$\\

 We define now a family of weighted shifts $\mathbb{H}_{p}$ acting on $\mathfrak{F}_{\beta}$ as following\\

$\mathbb{H}_{p} = \mathbb{M}^{p}\mathbb{D}^{p+1}$ with domain $D(\mathbb{H}_{p}) = \{\phi \in \mathfrak{F}_{\beta}; \mathbb{H}_{p}\phi \in \mathfrak{F}_{\beta}\}\hfill { } (2.20)$\\

\noindent Then we get\\

$\mathbb{H}_{p}^{*}e_{n}(z) = \mathbb{M}^{p+1}\mathbb{D}^{p}e_{n}(z) = \sqrt{m_{n+1}}\prod_{j=1}^{p}[m_{n-j+1}]e_{n+1}(z)$ for $ n \ge p \geq 0 $\\

\noindent i.e. $ \mathbb{H}_{p} ^{*}$ is weighted shift with weight $ \omega_{n,p} = \sqrt{m_{n+1}}\prod_{j=1}^{p}[m_{n-j+1}]$ for $ n = 1, .....$\\
and as we have denoted $[m_{n}]! = m_{1}.m_{2}.......m_{n} $ then $\omega_{n,p} = \sqrt{m_{n+1}}\frac{[m_{n}]!}{[m_{n-p}]!}$ for $ n \ge p \geq 0 $\\

\begin{Remark}
 (i) If $\beta \neq 2$ and p = 0  then the operator $\mathbb{H}_{0} = \mathbb{D} $ is particular case of Gelfond-Leontiev operator of generalized differentiation {\bf [10]} on $\mathfrak{F}_{\beta}$ and coincides with the usual differentiation on $\mathfrak{F}_{2}$.\\
(ii) For $\beta = 2$, It is known in {\bf[16]} that :\\
(a) the operator $\mathbb{H}_{p}$ with its domain $D(\mathbb{H}_{p} )$ is an operator chaotic on the classic Bargmann space.\\
(b) $ \mathbb{H}_{0} \phi_{\lambda}(z) = \lambda  \phi_{\lambda}(z) \quad \forall \quad \lambda \in I\!\!\!\!C$, where $ \phi_{\lambda}(z)  = \displaystyle {\sum_{n=0}^{\infty}}\frac{\lambda^{n}}{\sqrt{n!}} e_{n}(z)$ and\\
$\mid\mid\phi_{\lambda}\mid\mid^{2} = e^{\mid\lambda\mid^2}$\\
(c) The function $e^{-\mid\lambda\mid^2} \phi_{\lambda}(z) $ is called a coherent normalized quantum optics (see {\bf [26]})\\
(d) For $p = 1$, it is known that $ \mathbb{H}_{1} + \mathbb{H}_{1}^{*} $ is a not selfadjoint  operator and chaotic on the classic Bargmann space {\bf [8]}. This operator play an essential role in Reggeon field theory (see {\bf [20]}, {\bf [21]}and {\bf [22]})\\

\end{Remark}

D) \noindent Let $z = x +iy$, $\nu > 0$, $\Gamma = \mathbb{Z}\omega$
the discrete subgroup of the additive group $(I\!\!\!\!C, +)$ where $\omega \in \mathbb{C} -\{0\}$ and $\chi$ be a given map $\chi : \mathbb{Z}\omega \rightarrow U(1) =\{ \lambda \in I\!\!\!\!C ; \mid \lambda \mid = 1\}$ such that $\chi(\gamma) = e^{2i\pi\alpha m}$ for $\gamma = m\omega \in \mathbb{Z}\omega$, where $\alpha$ is a fixed real number. $\Gamma$ is non-compact lattice of the rank one and the triplet $(\nu, \Gamma, \chi )$ satisfies a Riemann-Dirac quantization type condition.\\
In {\bf [11]} Ghanmi and Intissar have considered the space $L_{\Gamma,\chi}^{2,\nu}(I\!\!\!\!C)$ of all measurable functions on $I\!\!\!\!C$ that are integrable on $\bigwedge(\Gamma)$ with respect $e^{-\nu\mid z\mid^{2}}dxdy$ and satisfying the function equation :\\

$\phi(z + \gamma)  = \chi(\gamma)e^{\nu z\overline{\gamma} + \frac{\nu\mid \gamma\mid^{2}}{2}}\phi(z)$ $\hfill { }(2.21)$\\

\noindent for almost every $z \in \mathbb{C}$ and every $\gamma \in \Gamma$,where $\bigwedge(\Gamma)$ is any given fundamental domain of the lattice $\Gamma$ witch is the strip $\mathbb{S} = [0, 1]\times\mathbb{R}$.\\

 $L_{\Gamma,\chi}^{2,\nu}(I\!\!\!\!C)$ is Hilbert space with the inner scalar product :\\

$< \phi_{1},  \phi_{2} >_{\Gamma} = \displaystyle {\int_{I\!\!\!\!C/\Gamma}}\displaystyle{\phi_{1}(z)\overline{\phi_{2}(z)}}e^{-\nu\mid z\mid^{2}}dxdy$ = $\displaystyle{{\int_{\bigwedge(\Gamma)}}\displaystyle{\phi_{1}(z)\overline{\phi_{2}(z)}}e^{-\nu\mid z\mid^{2}}dxdy}$  $\hfill { }(2.22)$\\

and associated norm:\\

$\displaystyle {\mid\mid \phi\mid\mid_{\Gamma} = \sqrt{\int_{\bigwedge(\Gamma)}\displaystyle{\mid\phi(z)\mid^{2}}e^{-\nu\mid z\mid^{2}}dxdy}}$ $\hfill { }(2.23)$\\

For $\omega = 1$, we write the function equation $(2.1)$ in the form:\\

$\phi(z + m)  = e^{2i\pi \alpha m}e^{\nu (z\gamma + \frac{m}{2})m}\phi(z)$ $\hfill { }(2.24)$\\

and we note $L_{\Gamma,\chi}^{2,\nu}(I\!\!\!\!C)$ by $L_{\Gamma,\alpha}^{2,\nu}(I\!\!\!\!C)$ \\

The $(\Gamma, \chi)$-theta Fock-Bargmann space $\mathfrak{F}_{\Gamma,\alpha}^{2,\nu}(I\!\!\!\!C)$ is defined now as a subspace of the space $\displaystyle{O(I\!\!\!\!C)}$ of holomorphic functions on $I\!\!\!\!C$, given by\\

$\displaystyle{\mathfrak{F}_{\Gamma,\alpha}^{2,\nu}(I\!\!\!\!C) = O(I\!\!\!\!C)\cap L_{\Gamma,\alpha}^{2,\nu}(I\!\!\!\!C) = \{\phi\in O(I\!\!\!\!C) ; < \phi, \phi >_{\Gamma} < \infty \}}$ $\hfill { }(2.25)$\\

Let \\

$\displaystyle{e_{n}^{\alpha, \nu}(z) = (\frac{2\nu}{\pi})^{1/4}e^{\frac{\nu}{2}z^{2}}e^{-\frac{\pi^{2}}{\nu}(n +\alpha)^{2} +2i\pi(n +\alpha)z}; n \in \mathbb{Z}}$ $\hfill { }(2.26)$\\

In theorem 2.5 of Ghanmi-Intissar {\bf [11]}, it is showed that the sequence $\displaystyle{e_{n}^{\alpha, \nu}(z) n \in \mathbb{Z}}$ is orthonormal basis of $\displaystyle{\mathfrak{F}_{\Gamma,\alpha}^{2,\nu}(I\!\!\!\!C)}$ and a function $\phi(z) = \displaystyle{\sum_{n \in \mathbb{Z}}a_{n}e_{n}^{\alpha, \nu}(z)}$ belongs to $\displaystyle{\mathfrak{F}_{\Gamma,\alpha}^{2,\nu}(I\!\!\!\!C)}$ if and only if $\displaystyle{\sum_{n \in \mathbb{Z}}\mid a_{n}\mid^{2} < +\infty}$\\

We consider now the operator of multiplication $\mathbb{M}$ by the function $M(z) = e^{2i\pi z}$ on the linear subspace $ \displaystyle{\mathfrak{F}_{\alpha}}$ of $\displaystyle{\mathfrak{F}_{\Gamma,\alpha}^{2,\nu}(I\!\!\!\!C)}$ generatedd by the orthogonal basis $\{e_{n}^{\alpha,\nu}(z); n \in\mathbb{N}\}$ (here we formallly take $e_{-1}^{\alpha,\nu}(z) = 0$) and $\mathbb{D}$ its adjoint.\\

The operator $\mathbb{M}$ acts on $e_{n}^{\alpha,\nu}(z)$; $n \in \mathbb{N}$ as following:\\

$\mathbb{M}e_{n}^{\alpha,\nu}(z) = \omega_{n}e_{n+1}^{\alpha,\nu}(z)$  $\hfill { } (2.27)$\\

with $\omega_{n} = c_{\alpha,\nu}e^{\frac{2\pi}{\nu}n}$ and $c_{\alpha, \nu}=e^{\frac{\pi}{\nu}+ 2\alpha}$\\

\noindent can be identified with\\

$\mathbb{M}(a_{n})_{n\in \mathbb{N}} = (\omega_{n}a_{n+1})_{n\in \mathbb{N}}$ $\hfill { } (2.28)$\\

\n and its adjoint is  given by :\\

$\mathbb{D}e_{n}^{\alpha,\nu}(z) = \omega_{n-1}e_{n-1}^{\alpha,\nu}(z)$  $\hfill { } (2.29)$\\

$\mathbb{D}(a_{n})_{n\in \mathbb{N}} = (\omega_{n-1}a_{n-1})_{n\in \mathbb{N}}$; $\omega_{-1} = 0$ $\hfill { } (2.30)$\\
\begin{Remark}
 {\it Let $p \in \mathbb{N}$ then $\mathbb{D}^{p+1}e_{p}^{\alpha,\nu}(z) = 0$.}\\
 \end{Remark}

 We define a family of unilateral weighted shifts $\mathbb{H}_{p}$ acting on $\displaystyle{\mathfrak{F}_{\alpha}}$ as following\\

$\mathbb{H}_{p} = \mathbb{M}^{p}\mathbb{D}^{p+1}$ with domain $D(\mathbb{H}_{p}) = \{\phi \in \displaystyle{\mathfrak{F}_{\alpha}}; \mathbb{H}_{p}\phi \in \displaystyle{\mathfrak{F}_{\alpha}}\}\hfill { } (2.31)$\\

\noindent Then we get\\

$\displaystyle{\mathbb{H}_{p}e_{n}^{\alpha,\nu}(z) = \mathbb{M}^{p}\mathbb{D}^{p+1}e_{n}^{\alpha,\nu}(z) = \omega_{n-1}[\prod_{j=1}^{p}\omega_{n-1-j}]^{2}e_{n-1}^{\alpha,\nu}(z)}$ $\hfill { } (2.32)$\\

and its adjoint is\\

$\displaystyle{\mathbb{H}_{p}^{*}e_{n}^{\alpha,\nu}(z) = \mathbb{M}^{p+1}\mathbb{D}^{p}e_{n}^{\alpha,\nu}(z) = \omega_{n}[\prod_{j=1}^{p}\omega_{n-j}]^{2}e_{n-1}^{\alpha,\nu}(z)}$ $\hfill { } (2.33)$\\

In following, if we put:\\

$\displaystyle{\mu := \frac{2\pi}{\nu}}$ and $\displaystyle{c_{\alpha} := c_{\alpha, \nu} =  e^{\frac{1}{2}(\mu + 4\alpha)}}$ and\\

$\displaystyle{\omega_{n,p} := \omega_{n}[\prod_{j=1}^{p}\omega_{n-j}]^{2}}$ $\hfill { } (2.34)$\\

we obtain\\

$\displaystyle{\omega_{n,p} = (c_{\alpha})^{2p +1}e^{(2p +1)\mu n}}$ $\hfill { } (2.35)$\\

and \\

$\displaystyle{\mathbb{H}_{p}e_{n}^{\alpha,\nu}(z) = \omega_{n-1,p}e_{n-1}^{\alpha,\nu}(z)}$ $\hfill { } (2.36)$\\

with $\displaystyle{\mathbb{H}_{p}e_{n}^{\alpha,\nu}(z) = 0}$ for $p \leq n$ and\\

$\displaystyle{\mathbb{H}_{p}^{*}e_{n}^{\alpha,\nu}(z) = \omega_{n,p}e_{n+1}^{\alpha,\nu}(z)}$ $\hfill { } (2.37)$\\

E) There exist many situations in physic where the Poincar\'e disk,\\ $\mathcal{D} = \{z \in \mathbb{C} ; |z| < 1\}$, is involved as a fundamental model or at least is used as a pedagogical toy (see for example Elwassouli et al in {\bf[9]}). It is a model of phase space for the motion of a material particle on a one sheeted two-dimensional hyperboloid viewed as a (1+1)-dimensional space-time with negative constant curvature, namely, the two dimensional anti de Sitter space-time.\\

The unit disk equipped with a K$\ddot{a}$hlerian potential,
$K_{\mathcal{D}}(z, \overline{z}) = \frac{1}{\pi}(1 -|z|^{2})^{2}$,
has the structure of a two-dimensional K$\ddot{a}$hlerian. Any K$\ddot{a}$hlerian manifold is
symplectic and so can be given a sense of phase space for some mechanical system.\\

Now let $\nu > \frac{1}{2}$ be a real parameter and let us equip the unit disk\\
 $\mathfrak{D} = \{ z \in \mathbb{C}; \mid z \mid < 1\}$ with a measure $d\lambda_{\nu}(z) = \frac{2\nu -1}{\pi}\frac{dxdy}{(1 -\mid z\mid^{2})^{2}}$.\\

For $\nu > 1$, we consider the Hilbert space $L_{\nu}^{2}(\mathfrak{D}, d\mu_{\nu}(z))$ of all functions $\phi$ on $\mathfrak{D}$ that are square integrable with respect to \\

$\displaystyle{d\mu_{\nu}(z) = (1 -\mid z\mid^{2})^{2\nu -2}dxdy}$ $\hfill { } (2.38)$\\

and introduce the Fock-Bargmann Hilbert space \\

$\displaystyle{\mathfrak{F}\mathfrak{B}_{\nu} = \mathcal{O}(\mathfrak{D})\cap L_{\nu}^{2}(\mathfrak{D}, d\mu_{\nu}(z))}$ $\hfill { } (2.39)$\\

where $\displaystyle{\mathcal{O}(\mathfrak{D})}$ is the space of all analytic functions $\phi(z)$ on $\mathfrak{D}$.

Let $z=re^{i\theta}$ with $0 < r < 1$ and $\theta\in [0, 2\pi]$, with respect to the scalar product defined on the holomorphic polynomials by:\\

$\displaystyle{<z^{n}, z^{m}> = \int_{\mathfrak{D}}z^{n}\overline{z^{m}}(1-\mid z\mid^{2})^{2\nu -2}dxdy} ; n \in \mathbb{N} , m \in \mathbb{N}$ $\hfill { } (2.40)$\\

we have:\\

$\displaystyle{\int_{\mathfrak{D}}z^{n}\overline{z^{m}}(1-\mid z\mid^{2})^{2\nu -2}dxdy = \int_{0}^{2\pi}\int_{0}^{1}r^{n+m+1}e^{i(n-m)\theta}drd\theta}$\\

$=\displaystyle{2\pi\int_{0}^{1}t^{n}(1- t)^{2\nu -2}dt}$\\

$ = \left\{\begin{array}[c]{l}\displaystyle{\frac{2\pi}{2\nu -1}\frac{\Gamma(2\nu)}{\Gamma(2\nu + n)}\Gamma(n+1)}\quad if\quad m=n \quad \\ \quad \\0 \quad  if \quad m \neq n\\
\end{array}\right.\hfill { }  (2.41)$\\

 This scalar product have the following property:\\

 the adjoint of the operator of differentiation \\

 $\displaystyle{\mathbb{A} = \frac{d}{dz}}$ is $\displaystyle{\mathbb{A}^{*} = z^{2}\frac{d}{dz} + 2\nu z}$ $\hfill { }(2.42)$\\

Then we associate to $\mathfrak{F}\mathfrak{B}_{\nu} $ the scalar product\\

$\displaystyle{< \phi, \psi>_{\mathfrak{D}} = \frac{2\nu -1}{2\pi} \int_{\mathfrak{D}}\phi(z)\overline{\psi(z)}d\mu_{\nu}(z)}$ $\hfill { } (2.43)$\\

to get\\

$\displaystyle{\mathbb{P}_{n}(z) = \sqrt{\frac{(2\nu)_{n}}{n!}}z^{n}; n \in \mathbb{N}}$ is an orthonormal basis of $\mathfrak{F}\mathfrak{B}_{\nu}$. $\hfill { } (2.44)$\\

where $\displaystyle{(2\nu)_{n} = \frac{\Gamma(2\nu + n)}{\Gamma(2\nu)}}$ is the Pochhammer symbol.\\

As on $\mathfrak{F}\mathfrak{B}_{\nu}$ the adjoint operator of \\

$\displaystyle{\mathbb{A} = \frac{d}{dz}}$ is the differential operator $\displaystyle{\mathbb{A}^{*} = z^{2}\frac{d}{dz} + 2\nu z}$\\

then they act on the orthonormal basis $\mathbb{P}_{n}(z)$ as following\\

$\displaystyle{\mathbb{A}^{*}\mathbb{P}_{n}(z) = \omega_{n}\mathbb{P}_{n+1}(z)}$ $\hfill { } (2.45)$\\

with $\displaystyle{\omega_{n} = \sqrt{(n+1)(2\nu + n)}}$.\\

and\\

$\mathbb{A}\mathbb{P}_{n}(z) = \omega_{n-1}\mathbb{P}_{n-1}(z)$ where $\mathbb{A}\mathbb{P}_{0}(z) = 0$ $\hfill { } (2.46)$\\

We define now a family of unilateral weighted shifts $\mathbb{H}_{p}$ acting on\\

$\displaystyle{\mathfrak{F}\mathfrak{B}_{\nu} = \mathcal{O}(\mathfrak{D})\cap L_{\nu}^{2}(\mathfrak{D}, d\mu_{\nu}(z))}; \nu > 1$ as following\\

$\mathbb{H}_{p} = \displaystyle{\mathbb{A}^{*^{p}}\mathbb{A}^{p+1}; p\in \mathbb{N}}$ be the linear unbounded densely defined shift operator acting on the $\displaystyle{\mathfrak{F}\mathfrak{B}_{\nu}}$ with domain \\

$D(\mathbb{H}_{p}) = \{\phi \in \displaystyle{\mathfrak{F}\mathfrak{B}_{\nu}} ; \mathbb{H}_{p}\phi \in \displaystyle{\mathfrak{F}\mathfrak{B}_{\nu}}\}$\\

whose its adjoint is defined by :\\

$\mathbb{H}^{*}_{p}\mathbb{P}_{n} = \omega_{n,p}\mathbb{P}_{n+1}$ $\hfill { } (2.47)$\\

where \\

$\displaystyle{\omega_{n,p} = \omega_{n}\prod_{j=1}^{p}\omega_{n-j}^{2}}; n \geq p\geq 0$ $\hfill { }  (2.48)$\\

i.e \\

$\displaystyle{\omega_{n,p} = \sqrt{(n+1)(2\nu + n)}\prod_{j=1}^{p} (n-j+1)(2\nu + n-j)}$ $\hfill { }  (2.49)$\\

\begin{Remark}

For the above spaces we will consider only \\

(i) the weights $\displaystyle{\omega_{n,p} = \sqrt{n+1}\prod_{j=1}^{p}(n-j+1)}$ for $ n \ge p \geq 0 $\\
associated to shifts acting on classic Bargmann space by noting that in this case\\

$\displaystyle{\omega_{n,p} \sim  n^{p +\frac{1}{2}}}$ $\hfill { }  (2.50)$\\

(ii) the weights $\displaystyle{\omega_{n,p} = \sqrt{m_{n+1}}\frac{[m_{n}]!}{[m_{n-p}]!}}$ for $ n \ge p \geq 0 $\\

where $m_{0} = 0$, $\displaystyle{[m_{n}]! = m_{1}.m_{2}.......m_{n}}$ and $\displaystyle{m_{n} = \frac{\Gamma(\frac{2}{\beta}(n+1))}{\Gamma(\frac{2}{\beta}n)}}$\\

associated to shifts acting on generalized Bargmann space by noting that in this case\\

$\displaystyle{\omega_{n,p} \sim  n^{\frac{2p +1}{\beta}}}$ $\hfill { }  (2.51)$\\

(iii) the weights $\displaystyle{\omega_{n,p} = (c_{\alpha})^{2p +1}e^{(2p +1)\mu n}}$ for $ n \ge p \geq 0 $ $\hfill { } (2.52)$\\

where\\

$\displaystyle{\mu := \frac{2\pi}{\nu}}$ and $\displaystyle{c_{\alpha} := c_{\alpha, \nu} =  e^{\frac{1}{2}(\mu + 4\alpha)}}$\\

associated to shifts acting on theta- Fock-Bargmann space.\\

(iv) the weights $\displaystyle{\omega_{n,p} = \sqrt{(n+1)(2\nu + n)}\prod_{j=1}^{p} (n-j+1)(2\nu + n-j)}$ for $ n \ge p \geq 0 $\\

associated to shifts acting on Fock-Bargmann space on Poincar\'e disk by noting that in this case\\

$\displaystyle{\omega_{n,p} \sim n^{2p +1}}$ $\hfill { }  (2.53)$\\

\end{Remark}

\section{ On chaoticity of the sum of chaotic shift and its adjoint and applications}

In this section we give sufficient conditions on a linear unbounded densely defined chaotic shift operator $\mathbb{T}$ acting on a Hilbert space such that $\mathbb{T} + \mathbb{T}^{*}$ is chaotic where $\mathbb{T}^{*}$ is its adjoint.\\

\begin{Theorem}

Let a linear unbounded densely defined chaotic shift operator $(\mathbb{T}, D(\mathbb{T}))$ on a Hilbert space $\mathbb{E} = \displaystyle{\{\phi; \phi = \sum_{n=1}^{\infty}a_{n}e_{n}\}}$ such that its adjoint is defined by:\\

$\mathbb{T}^{*}e_{n} = \omega_{n}e_{n+1}$ $\hfill { } (3.1)$\\

where $\{e_{n}\}$ is an orthonormal basis of $\mathbb{E}$ and $\omega_{n}$ is positive weight associated to $\mathbb{T}$\\

We assume that\\

(Assumption $Hyp_{1}$)\quad $\displaystyle{\sum_{n=1}^{\infty}\frac{1}{\omega_{n}} < \infty}$ $\hfill { } (3.2)$\\

(Assumption $Hyp_{2}$)\quad $\displaystyle{\omega_{n-1}\omega_{n+1} \leq \omega_{n}^{2}}$ $\hfill { } (3.3)$\\

(Assumption $Hyp_{3})$\quad there exist $\alpha > 0$, $\beta > 0$, $a > 0$, and a sequence $\gamma_{n}$ that:\\

(1) $\displaystyle{\frac{\omega_{n}\gamma_{n}}{\gamma_{n+1}} \geq n^{1+\alpha}}$  $\hfill { } (3.4)$\\

(2) $\displaystyle{\frac{\omega_{n-1}\gamma_{n+1}}{\omega_{n}\gamma_{n-1}} = 1 - \frac{a}{n} + O(\frac{1}{n^{1+\beta}})}$ $\hfill { } (3.5)$\\

and\\

(3) $\displaystyle{\sum_{k=1}^{\infty}\frac{1}{\gamma_{n}^{2}} < \infty}$ $\hfill { } (3.6)$\\

Then  for $\lambda \in \mathbb{C}$ the following recurrence sequence\\

$ (*) \left\{\begin{array}[c]{l}\displaystyle{u_{1}(\lambda) = 1}\\
\displaystyle{u_{2}(\lambda) = \frac{\lambda}{\omega_{1}}} \\
\displaystyle{\omega_{n-1}u_{n-1}(\lambda) + \omega_{n}u_{n+1}(\lambda) = \lambda u_{n}(\lambda)}\\
\end{array}\right.\hfill { }  (3.7)$\\

\quad\\

(i) is solvable  for all $\lambda \in \mathbb{C}$.\\

(ii) $\displaystyle{\sum_{n=1}^{\infty}\mid u_{n}(\lambda)\mid^{2} < \infty}$ for all $\lambda \in \mathbb{C}$.\\

(iii) the spectrum of $\mathbb{T} + \mathbb{T}^{*}$ is the all complex plane $\mathbb{C}$.\\

(iv) $(\mathbb{T} + \mathbb{T}^{*})^{m}$ is closed $\quad \forall \quad m \in \mathbb{N}$.\\

(v) $\mathbb{T} + \mathbb{T}^{*}$ is hypercyclic operator.\\

(vi) $\mathbb{T} + \mathbb{T}^{*}$ is chaotic operator.\\
\end{Theorem}

{\bf Proof of theorem}\\

(i) By using the Yu. Berzanskii'theory on the difference operators in {\bf[4]} in particular the theorem $1.5$ ch. $VII$ and the above asymptions (3.2) and (3.3) we deduce that the sequence $u_{n}(\lambda)$ is always solvable and is a polynomial of degree $n-1$ called the polynomials of first kind associated to the operato $\mathbb{T}$.\\

(ii) Let $M > 0$ (large enough) such that $ \displaystyle{\mid u_{n}(\lambda)\mid \leq \frac{M}{\gamma_{n}}}$ and $ \displaystyle{\mid u_{n-1}(\lambda)\mid \leq \frac{M}{\gamma_{n-1}}}$.

As $\displaystyle{\omega_{n-1}u_{n-1}(\lambda) + \omega_{n}u_{n+1}(\lambda) = \lambda u_{n}(\lambda)}$ Then\\

$\displaystyle{u_{n+1}(\lambda) = \frac{\lambda }{\omega_{n}}u_{n}(\lambda) - \frac{\omega_{n-1}}{\omega_{n}}u_{n-1}(\lambda)}$ then \\

$\displaystyle{\mid u_{n+1}(\lambda)\mid \leq \frac{\mid\lambda\mid}{\omega_{n}}\mid u_{n}(\lambda)\mid + \frac{\omega_{n-1}}{\omega_{n}}\mid u_{n-1}(\lambda)\mid}$\\

$\displaystyle{\leq M[\frac{\mid\lambda \mid}{\omega_{n}}\frac{1}{\gamma_{n}} + \frac{\omega_{n-1}}{\omega_{n}}\frac{1}{\gamma_{n-1}}]}$\\

$\displaystyle{\leq \frac{M}{\gamma_{n+1}}[\frac{\mid\lambda \mid}{\omega_{n}}\frac{\gamma_{n+1}}{\gamma_{n}} + \frac{\omega_{n-1}}{\omega_{n}}\frac{\gamma_{n+1}}{\gamma_{n-1}}]}$\\

From (3.4) and (3.5) we get\\

$\displaystyle{\mid u_{n+1}(\lambda)\mid \leq \frac{M}{\gamma_{n+1}}[\frac{\mid \lambda\mid}{n^{1+\alpha}} + 1 - \frac{a}{n} + O(\frac{1}{n^{1+\beta}})]}$\\

$\displaystyle{\mid u_{n+1}(\lambda)\mid \leq \frac{M}{\gamma_{n+1}}[ 1 - \frac{a}{n} +\frac{\mid \lambda\mid}{n^{1+\alpha}} + O(\frac{1}{n^{1+\beta}})]}$\\

$\displaystyle{\leq \frac{M}{\gamma_{n+1}}}$\\

and from (3.6) we deduce that \\

$\displaystyle{\sum_{n=1}^{\infty}\mid u_{n}(\lambda)\mid^{2} < \infty}$ for all $\lambda \in \mathbb{C}$\\

(iii) Let $u_{n}(\lambda)$ the sequence defined by (3.7) and $\phi_{\lambda} = \displaystyle{\sum_{n=1}^{\infty} u_{n}(\lambda)e_{n}}$ then $(\mathbb{T} + \mathbb{T}^{*})\phi_{\lambda} = \lambda \phi_{\lambda}$.\\

 As $\displaystyle{\sum_{n=1}^{\infty}\mid u_{n}(\lambda)\mid^{2} < \infty}$ for all $\lambda \in \mathbb{C}$ then the spectrum of $\mathbb{T} + \mathbb{T}^{*}$ is $\mathbb{C}$\\

(iv) As $\mathbb{T}$ is chaotic we have $\mathbb{T}^{m}$, $\mathbb{T}^{*^{m}}$ and $\mathbb{T}^{k}\mathbb{T}^{*^{j}}$ are closed operators $\quad \forall\quad (m,k,j) \in \mathbb{N}^{3}$ then $(\mathbb{T} + \mathbb{T}^{*})^{{m}}$ is closed.\\

(v) We verify now that the operator $\mathbb{T} + \mathbb{T}^{*}$ on $\mathbb{E}$ satisfies the hypercyclicity criterion, as quoted above. \\

Let $\Omega_{1} = \{\lambda \in \mathbb{C} ; \mid \lambda \mid > 1\}$, $\Omega_{2} = \{\lambda \in \mathbb{C} ; \mid \lambda \mid < 1\}$, $\Omega_{3} = \{\lambda \in \mathbb{C} ; \mid \lambda \mid = 1\}$\\

and\\

$\mathbb{F}_{1}$ the space spanned by $ \displaystyle{\{\phi_{\lambda} = \sum_{n=1}^{\infty}u_{n}(\lambda)e_{n} ; \lambda \in \Omega_{1}\}}$\\

$\mathbb{F}_{2}$ the space spanned by $ \displaystyle{\{\phi_{\lambda} = \sum_{n=1}^{\infty}u_{n}(\lambda)e_{n} ; \lambda \in \Omega_{2}\}}$\\

$\mathbb{F}_{3}$ the space spanned by $ \displaystyle{\{\phi_{\lambda} = \sum_{n=1}^{\infty}u_{n}(\lambda)e_{n} ; \lambda \in \Omega_{3}\}}$\\

each $\mathbb{F}_{j}$ $j=1, 2,3$ is dense in $\mathbb{E}$ because if $\phi \in \mathbb{F}_{j}^{\perp}$ then the function $g(\lambda) = <\phi, \phi_{\lambda}>$ is entire on $\mathbb{C}$ and equal to zero on $\mathbb{F}_{j}$ witch has accumulation points, then $g(\lambda) = 0$ on $\mathbb{C}$ and also $\phi = 0$.\\

Let $\mathbb{S}$ be the linear mapping on $\mathbb{F}_{1}$ determined by:\\

$\displaystyle{\mathbb{S}\phi_{\lambda} = \mathbb{S}(\sum_{n=1}^{N}a_{n}\phi_{\lambda_{n}})  = \sum_{n=1}^{N}\frac{a_{n}}{\delta_{n}}\phi_{\lambda_{n}}}; \lambda_{n} \in \omega_{1}$ and $a_{n} \in \mathbb{C}$\\

then $(\mathbb{T} + \mathbb{T}^{*})^{n} \rightarrow 0$ pointwise on $\mathbb{F}_{1}$, and $\mathbb{T}\mathbb{S} =I$ and $S^{n} \rightarrow 0$ pointwise on $\mathbb{F}_{2}$. So with the property (iv) we deduce that $\mathbb{T} + \mathbb{T}^{*}$ is hypercyclic on $\mathbb{E}$ by using the lemme of Godefroy-Shapiro or the theorem of B$\grave{e}$ et al.\\

(vi) For see that $\mathbb{F}_{3}$ is subset of periodic points of $\mathbb{T} + \mathbb{T}^{*}$, we take $a_{m} \in \mathbb{C}$, $n_{m},k_{m} \in \mathbb{Z}$, $\delta_{m} = \displaystyle{e^{2i\pi\frac{ n_{m}}{k_{m}}}}$ and $\displaystyle{\phi = \sum_{m=1}^{N}a_{m}\phi_{\delta_{m}}}$. Then we observe that for $l = \prod_{m=1}^{N}k_{m}$ we have $(\mathbb{T} + \mathbb{T}^{*})^{l}\phi = \phi$.\\

\begin{Theorem}
Let $\displaystyle{\mathbb{B} = \{\phi : \mathbb{C} \rightarrow \mathbb{C} \quad entire \quad ; \int_{\mathbb{C}}\mid \phi(z)\mid^{2}e^{-\mid z\mid^{2}}dxdy < \infty\}}$ the classic Bargmann space and $\mathbb{H}_{p} = \displaystyle{z^{p}\frac{d^{p+1}}{dz^{p+1}}}; p\in \mathbb{N}$ the linear unbounded densely defined shift operator acting on the $\mathbb{B}_{p}$ space; $p=0, 1, ...$\\

 $\displaystyle{\mathbb{B}_{p}=\{\phi \in \mathbb{B}; \frac{d^{j}}{dz^{j}}\phi(0)= 0; 0 \leq j \leq p \}}$\\

with domain \\

$D(\mathbb{H}_{p}) = \{\phi \in B ; \mathbb{H}_{p}\phi \in B\}\cap B_{p}$\\

whose its adjoint is defined by :\\

$\mathbb{H}^{*}_{p}e_{n} = \omega_{n,p}e_{n+1}$ $\hfill { } ((3.8)$\\

where \\

$ \displaystyle{\omega_{n,p}} = \left\{\begin{array}[c]{l}\displaystyle{\sqrt{n+1}}\quad\quad \quad\quad \quad for\quad p = 0\\
\displaystyle{\sqrt{n+1}\prod_{j=0}^{p-1}(n-j)} \quad for \quad p\geq 1 \\
\end{array}\right.\hfill { }  (3.9)$\\

\quad\\

and\\

$\displaystyle{\{e_{n}(z) = \frac{z^{n}}{\sqrt{n!}}; n = p, p+1, ....\}}$ is an orthonormal basis of $\mathbb{B}_{p}$\\

Then we have\\

(i) for all $p \geq 0, \mathbb{H}_{p}$ is chaotic.\\

(ii) for all $p \geq 1, \mathbb{H}_{p} + \mathbb{H}^{*}_{p}$ is chaotic.\\

\end{Theorem}

{\bf Proof}\\

(i) In {\bf[16]}, the author showed that $\mathbb{H}_{p}$ acting on $\mathbb{B}_{p}$ is chaotic in the sense of Devaney for all $p \in \mathbb{N}$.\\

(ii) We will be concerned with the chaoticity of $\mathbb{H}_{p} + \mathbb{H}_{p}^{*}$ for $p \geq 1$ by using the above theorem. We begin to remark that\\

-the assumption  $(3.2)$ of the above theorem is verified because
$\displaystyle{\omega_{n,p}\sim n^{p+1/2}}$ and as $p \geq 1$ we have \\

$\displaystyle{\sum_{n=1}^{\infty}\frac{1}{\omega_{n,p}} < \infty}$.\\

-to verify the assumption $(3.3)$ of the above theorem, we write $\omega_{n,p}$ in the following form \\

 $\displaystyle{\omega_{n,p} = \sqrt{n+1}n(n-p+1)A_{n,p}}$ where $\displaystyle{A_{n,p} = \prod_{j=1}^{p-2}(n-j)}$\\

 then\\

 $\displaystyle{\omega_{n-1,p} = \sqrt{n}(n-p+1)(n-p)A_{n,p}}$\\

 and\\

 $\displaystyle{\omega_{n+1,p} = \sqrt{n+1}n(n+1)A_{n,p}}$\\

 and as $n \geq p$ we deduce that\\

- $\displaystyle{\omega_{n-1,p}\omega_{n+1,p}\leq \omega_{n-1,p}^{2}}$ for all $n \geq p$\\

- $\displaystyle{\frac{\omega_{n-1,p}}{\omega_{n,p}} = \frac{\sqrt{n}(n-p)}{\sqrt{n+1}n}}$\\

We choose choose\\

$\gamma_{n} = \sqrt{n}Logn$ $\hfill { } (3.10)$\\

to obtain :\\

- $\displaystyle{\frac{\mid\lambda \mid}{\omega_{n,p}}\frac{\gamma_{n+1}}{\gamma_{n}} = \frac{\mid\lambda \mid}{\omega_{n,p}}\frac{\sqrt{n+1}Log(n+1)}{\sqrt{n}Logn} \sim \frac{1}{n^{p+ 1/2}}}$\\

- $\displaystyle{\frac{\omega_{n-1,p}}{\omega_{n,p}}\frac{\gamma_{n+1}}{\gamma_{n-1}} = \frac{n-p}{\sqrt{n}\sqrt{n-1}}\frac{Log(n+1)}{Log(n-1)}}$\\

$= \displaystyle{(1-\frac{p}{n})(1 -\frac{1}{n})^{-\frac{1}{2}}\frac{Log(n+1)}{Log(n-1)}}$\\

As $\displaystyle{(1-\frac{p}{n})(1 -\frac{1}{n})^{-\frac{1}{2}} \leq 1 - \frac{2p -1}{2n} + O(\frac{1}{n^{2}} )}$\\

 and\\

 $\displaystyle{ \frac{Log(n+1)}{Log(n-1)} = \frac{Log(n-1+2)}{Log(n-1)} = \frac{Log(n-1) + Log(1 + \frac{2}{n-1})}{Log(n-1)}}$ \\

 then there exist $ C > 0; \displaystyle{ \frac{Log(n+1)}{Log(n-1)} \leq 1 + \frac{C}{(n-1)Log(n-1)}}$ \\

 we get there exist $a > 0$ and $\beta > 0$ such that\\

$\displaystyle{\frac{\mid\lambda \mid}{\omega_{n,p}}\frac{\gamma_{n+1}}{\gamma_{n}} + \frac{\omega_{n-1,p}}{\omega_{n,p}}\frac{\gamma_{n+1}}{\gamma_{n-1}} \leq 1 - \frac{a}{n} + 0(\frac{1}{n^{1 + \beta}}) \leq 1}$ as $ n$ enough large.\\

this implies that $\mathbb{H}_{p} + \mathbb{H}_{p}^{*}$ is chaotic.\\

\begin{Remark}

(i) The operator $\displaystyle{\mu z^{p}\frac{d^{p}}{dz^{p}} + i\lambda z^{p}(\frac{d^{p}}{dz^{p}} + z^{p})\frac{d^{p}}{dz^{p}}}$ with $\mu > 0$, $\lambda \in \mathbb{R}$ and $i = \sqrt{-1}$ is not chaotic on $\mathbb{B}_{p}$, in fact it is an operator with compact resolvent.\\
For $p=1$, the Reggeon field theory is governed by this non-self adjoint operator see {\bf [21]} and references therein .\\

(ii) We have assumed the assumptions Hyp1 and Hyp2 in theorem 3.1 for the solvability of  the sequence defined by (3.7):\\

$ (*) \left\{\begin{array}[c]{l}\displaystyle{u_{1}(\lambda) = 1}\\
\displaystyle{u_{2}(\lambda) = \frac{\lambda}{\omega_{1}}} \\
\displaystyle{\omega_{n-1}u_{n-1}(\lambda) + \omega_{n}u_{n+1}(\lambda) = \lambda u_{n}(\lambda)}\\
\end{array}\right.$\\

\quad\\

(iii) the assumption Hyp3 in theorem 3.1 can be replaced by\\

There exist a sequence $\gamma_{n}$ that:\\

- $\displaystyle{\frac{\mid\lambda\mid}{\omega_{n}}\frac{\gamma_{n+1}}{\gamma_{n}} + \frac{\omega_{n-1}}{\omega_{n}}\frac{\gamma_{n+1}}{\gamma_{n-1}} < 1}$ ; $\lambda \in \mathbb{C}$ and $n$ large enough $\hfill { } (3.11)$\\

and\\

- $\displaystyle{\sum_{k=1}^{\infty}\frac{1}{\gamma_{n}^{2}} < \infty}$ $\hfill { } (3.12)$\\

\end{Remark}

As application of theorem 3.1 where the assumption Hyp3 is replaced by (3.11) and (3.12), we give next application\\

\begin{Theorem}
Let $\displaystyle{\mathfrak{F}_{\alpha}}$ the lattice Fock-Bargmann space and $\mathbb{H}_{p} = \displaystyle{\mathbb{M}^{p}\mathbb{D}^{p+1}}; p\in \mathbb{N}$ the linear unbounded densely defined shift operator acting on $\mathfrak{F}_{\alpha,p}$; $p=0, 1, ...$ the space spanned by $\displaystyle{e_{n}^{\alpha, \nu}(z) = (\frac{2\nu}{\pi})^{1/4}e^{\frac{\nu}{2}z^{2}}e^{-\frac{\pi^{2}}{\nu}(n +\alpha)^{2} +2i\pi(n +\alpha)z}; n = p, p+1, ....\}}$\\

with domain \\

$D(\mathbb{H}_{p}) = \{\phi \in \mathfrak{F}_{\alpha} ; \mathbb{H}_{p}\phi \in\mathfrak{F}_{\alpha}\}\cap \mathfrak{F}_{\alpha,p}$\\

whose its adjoint is defined by :\\

$\mathbb{H}^{*}_{p}e_{n}^{\alpha,\nu} = \omega_{n,p}e_{n+1}^{\alpha,\nu}$ $\hfill { } (3.13)$\\

where \\

$\displaystyle{\omega_{n,p} = (c_{\alpha})^{(2p +1)}e^{(2p+1)\mu n}}$ $\hfill { } (3.14)$\\

with $\displaystyle{\mu := \frac{2\pi}{\nu}}$ and $\displaystyle{c_{\alpha} =  e^{\frac{1}{2}(\mu + 4\alpha)}}$ \\

and\\

$\displaystyle{\{e_{n}^{\alpha,\nu}(z) =  (\frac{2\nu}{\pi})^{1/4}e^{\frac{\nu}{2}z^{2}}e^{-\frac{\pi^{2}}{\nu}(n +\alpha)^{2} +2i\pi(n +\alpha)z}; n = p, p+1, ....\}}$ is an orthonormal basis of $\mathfrak{F}_{\alpha,p}$\\

Then we have\\

(i) for all $p \geq 0, \mathbb{H}_{p}$ is chaotic.\\

(ii) for all $p \geq 0, \mathbb{H}_{p} + \mathbb{H}^{*}_{p}$ is chaotic.\\

\end{Theorem}

{\bf Proof}\\

(i) In {[18]}, we have shown that for all $p \geq 0, \mathbb{H}_{p}$ is chaotic.\\

(ii) Let $\displaystyle{\omega_{n,p} = (c_{\alpha})^{(2p +1)}e^{(2p+1)\mu n}}$,if we choose now\\

$\displaystyle{\gamma_{n} = (c_{\alpha})^{(2p +1)}e^{\frac{(2p+1)\mu}{\beta} n}}$ with $\beta > 2$ then we get\\

$\displaystyle{\frac{\mid \lambda\mid}{(c_{\alpha})^{(2p +1)}e^{(2p+1)\mu n}}\frac{\gamma_{n+1}}{\gamma_{n}} + \frac{(c_{\alpha})^{(2p +1)}e^{(2p+1)\mu (n-1)}}{(c_{\alpha})^{(2p +1)}e^{(2p+1)\mu n}}\frac{\gamma_{n+1}}{\gamma_{n-1}}} =$

$\displaystyle{\frac{\mid \lambda\mid}{(c_{\alpha})^{(2p +1)}e^{(2p+1)\mu n}}\frac{\gamma_{n+1}}{\gamma_{n}} + e^{-(2p +1)\mu}\frac{\gamma_{n+1}}{\gamma_{n-1}}} =$\\

$\displaystyle{\frac{\mid \lambda\mid}{(c_{\alpha})^{(2p +1)}e^{(2p+1)\mu n}}e^{(2p+1)\frac{\mu}{\beta}} + e^{-(2p +1)\mu (1 - \frac{\beta}{2})}}$\\

Let $m_{\beta} = \displaystyle{\frac{\mid \lambda\mid}{(c_{\alpha})^{(2p +1)}}e^{(2p+1)\frac{\mu}{\beta}}}$\\

 and \\

 $C_{\beta} = \displaystyle{e^{-(2p +1)\mu(1-\frac{\beta}{2})}}$\\

 then\\

For $\beta > 2$ and $n \geq \frac{1}{(2p + 1)\mu}Log\frac{m_{\beta}}{1 -C_{\beta}}$\\

 we deduce that \\

$\displaystyle{m_{\beta}e^{-(2p+1)\mu n} + C_{\beta} < 1}$ $\hfill { } (3.15)$\\

and\\

$\displaystyle{\sum_{n=1}^{\infty}\mid u_{n}(\lambda)\mid^{2} < \infty}$ $\hfill { } (3.16)$\\

because $\displaystyle{\sum_{n=1}^{\infty}\frac{1}{\gamma_{n,p}^{2}} < \infty}$\\

the verification of the rest of hypothesis of theorem 3.1 for $\omega_{n,p}$ are obvious.\\ This shows that $\mathbb{H}_{p} + \mathbb{H}_{p}^{*}$ is chaotic for $p=0, 1, 2, ......$\\

\begin{Remark}

(i) In each of the above cases the expansion coefficients $u_{n}(\lambda)$ satisfy a three-term recurrence relation of the form\\

$u_{n+1} + \alpha_{n}u_{n} + \beta_{n}u_{n-1} = 0$ $\hfill{ } (*)$\\

 We can consider the above equation as a second order linear homogeneous
difference equation.\\

let us assume that $lim \quad\alpha_{n} = \alpha$ and $lim \quad\beta_{n} =\beta$ as $n \rightarrow \infty$. Then the Poincar\'e analysis holds see {\bf [28]} or {\bf [29]} and states that if $t_{l}$ and $t_{2}$ are the two roots of the quadratic equation
$t^{2} + \alpha t + \beta = 0$ and $\mid t_{2} \mid > \mid t_{1}\mid $ then (*) possesses two linearly independent solutions \\

This Poincar\'e analysis is not applicable to \\

$\displaystyle{u_{n+1}(\lambda) - \frac{\lambda}{\omega_{n,p}}u_{n}(\lambda) +  \frac{\omega_{n-1,p}}{\omega_{n,p}} u_{n-1}(\lambda) = 0}$ $\hfill{ } (**)$\\

 we have $\displaystyle{lim \frac{\lambda}{\omega_{n,p}} = \alpha = 0}$ and $\displaystyle{lim \frac{\omega_{n-1,p}}{\omega_{n,p}} = \beta = 1}$ as $n \rightarrow \infty$.\\

but the two roots $ t_{2} = i$ and $ t_{1} = -i $ of the quadratic equation $t^{2} + \alpha t + \beta = 0$ have same modulus.\\

(ii) In {\bf[23]}, we study directly the operator $\mathbb{H}_{0} + \mathbb{H}_{0}^{*}$ on classic Bargmann space and we give the application of theorem 3.1 for $\mathbb{H}_{p} + \mathbb{H}_{p}^{*}$ on generalized Bargmann space with $\beta \neq 2$. This last application use the following fondamental lemma.
\end{Remark}

\begin{Lemma}
{\it Let $ m_{0} = 0, m_{n} = \frac{\Gamma(\frac{2}{\beta}(n+1))}{\Gamma(\frac{2n}{\beta})}$;\quad $ n =1,2, ...$ then \\

$m_{n}\sim (\frac{2}{\beta})^{\frac{2}{\beta}}n^{\frac{2}{\beta}}, \quad n \rightarrow +\infty$} $\hfill { } (3.17)$\\
\end{Lemma}

\section{Appendix: Chaoticity of $\mathbb{H}_{p}$ on Fock-Bargmann space associated to Poincar\'e disk}

In this appendix, we show  that the operators $\mathbb{H}_{p} = \mathbb{A}^{*^{p}}\mathbb{A}^{p+1}; p\in \mathbb{N}$ with  domain $D(\mathbb{H}_{p}) = \{\phi \in \mathfrak{F}\mathfrak{B}_{\nu}; \mathbb{H}_{p}\phi \in \mathfrak{F}\mathfrak{B}_{\nu}\}$ are chaotic, where $\mathfrak{F}\mathfrak{B}_{\nu}$ is Fock-Bargmann space associated to Poincar\'e disk, $\displaystyle{\mathbb{A} = \frac{d}{dz}}$ and $\displaystyle{\mathbb{A}^{*} = z^{2}\frac{d}{dz} + 2\nu z}$ its adjoint.\\

\begin{Theorem}
{\it Let $\displaystyle{\mathfrak{F}\mathfrak{B}_{\nu} = \mathcal{O}(\mathfrak{D})\cap L_{\nu}^{2}(\mathfrak{D}, d\mu_{\nu}(z))}; \nu > 1$ }\\

{\it $\mathbb{H}_{p} = \displaystyle{\mathbb{A}^{*^{p}}\mathbb{A}^{p+1}; p\in \mathbb{N}}$ be the linear unbounded densely defined shift operator acting on the $\mathfrak{F}\mathfrak{B}_{\nu,p}$ the space spanned by orthonormal basis} \\

$\displaystyle{\mathbb{P}_{n}(z) = \sqrt{\frac{(2\nu)_{n}}{n!}}z^{n}; n = p, p+1, ....\}}$$\hfill { } (4.1)$\\

{\it where $\displaystyle{(2\nu)_{n} = \frac{\Gamma(2\nu + n)}{\Gamma(2\nu)}}$ is the Pochhammer symbol.}\\

{\it with domain }\\

$D(\mathbb{H}_{p}) = \{\phi \in \displaystyle{\mathfrak{F}\mathfrak{B}_{\nu}} ; \mathbb{H}_{p}\phi \in \displaystyle{\mathfrak{F}\mathfrak{B}_{\nu}}\}\cap \displaystyle{\mathfrak{F}\mathfrak{B}_{\nu,p}}$\\

{\it whose its adjoint is defined by :}\\

$\displaystyle{\mathbb{H}^{*}_{p}\mathbb{P}_{n} = \omega_{n,p}\mathbb{P}_{n+1}}$ $\hfill { } (4.2)$\\

{\it where }\\

$ \displaystyle{\omega_{n,p} = \sqrt{(n+1)(2\nu + n)}\prod_{j=1}^{p}[(n+1)(2\nu + n)]^{2}}; n \geq p \geq 0$ $\hfill { }  (4.3)$\\

{\it Then we have}\\

{\it For all $p \geq 0, \mathbb{H}_{p}$ is chaotic.}\\

\end{Theorem}

{\bf Proof}

To use the theorem of $B\grave{e}s$ and al, we begin by observing that for $\displaystyle{\phi(z) =\sum_{k=p}^{\infty}a_{k}\mathbb{P}_{k}(z)}$ such that $\displaystyle{\sum_{k=p}^{\infty}\mid a_{k}\mid^{2} < \infty}$ we have the obvious properties\\

 (i) $\displaystyle{\mathbb{H}_{p}^{m}\phi(z) =  \sum_{k=p}^{\infty}[\prod_{j=p}^{m+k-1}\omega_{j,p}]a_{k+m}\mathbb{P}_{k}(z)}$ of domain\\

$\displaystyle{D(\mathbb{H}_{p}^{m}) = \{ \phi = \sum_{k=p}^{\infty}a_{k}\mathbb{P}_{k};\sum_{k=p}^{\infty}\mid a_{k}\mid^{2} < \infty \quad and \quad \sum_{k=p}^{\infty}[\prod_{j=p}^{m+k-1}\omega_{j,p}]^{2}\mid a_{k+m}\mid^{2} < +\infty\}}$\\

 witch is dense in $\displaystyle{\mathfrak{F}\mathfrak{B}_{\nu,p}} \quad\forall\quad m \in \mathbb{N}$\\

(ii) $\mathbb{H}_{p}^{m}$ is closed $\quad\forall\quad m \in \mathbb{N}$ and $\mathbb{H}_{p}^{m} \mathbb{P}_{k}(z) = 0$ $\quad\forall\quad m  > k \geq p \geq 0$\\

(iii) As $\omega_{n,p} \rightarrow +\infty$ then the spectrum of $\mathbb{H}_{p}$ is the all complex plane.\\

In fact, let $\displaystyle{\phi_{\lambda} = \sum_{k=p}^{\infty}a_{k}\mathbb{P}_{k}(z)}$ with $\displaystyle{a_{k} = \prod_{j=p}^{k-1}\frac{\lambda}{\omega_{j,p}}}$ i.e $\displaystyle{\phi_{\lambda} =\sum_{k=p}^{+\infty}[\prod_{j=p}^{k-1}\frac{\lambda}{\omega_{j,p}}]P_{k}(z)}$ then as $a_{p} = 0$ we deduce that\\

$\mathbb{H}_{p}\phi_{\lambda} = \lambda\phi_{\lambda} , \forall \quad \lambda \in \mathbb{C}$.$\hfill { }  (4.4)$\\\\

and as
$\displaystyle{\sum_{k=p}^{\infty}[\prod_{j=p}^{k}\frac{\lambda}{\omega_{j,p}}]^{2} < +\infty}$
then $\phi_{\lambda} \in D(\mathbb{H}_{p})$\\

Now, take $\mathbb{Y}$ the linear subspace generated by finite combinations of basis $\{\mathbb{P}_{k}\}_{k=p}^{\infty}$, this subspace $\mathbb{Y}$ is dense in $\displaystyle{\mathfrak{F}\mathfrak{B}_{\nu,p}}$ and we define on it the operator $\mathbb{S}$ acting on $\displaystyle{\phi = \sum_{k=p}^{N}a_{k}\mathbb{P}_{k}}$ as following \\

$\displaystyle{\mathbb{S}\phi = \sum_{k=p}^{N+1}\frac{a_{k-1}}{\omega_{k-1,p}}\mathbb{P}_{k}}$ $\hfill { }  (4.5)$\\ then \\

$\displaystyle{\mathbb{S}^{n}\mathbb{P}_{k} = \frac{1}{\prod_{j=k}^{n+k}\omega_{j,p}}\mathbb{P}_{k+n}}$\\

as $\displaystyle{\prod_{j=p}^{n}\omega_{j,p} \rightarrow +\infty}$ as $n \rightarrow +\infty$ we get\\

$\displaystyle{\mathbb{S}^{n}\mathbb{P}_{k} \rightarrow 0}$ in $\displaystyle{\mathfrak{F}\mathfrak{B}_{\nu,p}}$ as $n \rightarrow +\infty$ $\hfill { }  (4.6)$\\\\

By noting that $\mathbb{H}_{p}^{n}\mathbb{P}_{k} = 0 $ for $n > k$ and any element of $\mathbb{Y}$ can be annihilated by a finite power of $\mathbb{H}_{p}$ and $\mathbb{H}_{p}\mathbb{S}_{p} = \mathbb{I}_{\mid \mathbb{Y}}$ then the hyperciclycity of $\mathbb{H}_{p}$ follows from the theorem of $ B\grave{e}s$ and al. recalled above.\\

We shall now show that $\mathbb{H}_{p}$ has a dense set of periodic points. To see this, it suffices to show that for every element $\phi$ in the dense subspace $Y$ there is a periodic point $\psi$ arbitrarily close to it.\\

For $s \geq p$ and $N \geq s$ we put\\

$\displaystyle{\varphi_{s,N}(z) = \mathbb{P}_{s}(z) + \sum_{k=s+1}^{\infty}[\prod_{j=s}^{kN + s -1}\frac{1}{\omega_{j,p}}]\mathbb{P}_{kN +s}(z)}$ $\hfill { }  (4.7)$\\

Then we have the following obvious lemma\\

\begin{Lemma}

(i) $\displaystyle{\mathbb{H}_{p}^{N}\prod_{j=0}^{kN -1}\frac{1}{\omega_{j,p}} \mathbb{P}_{kN} =\prod_{j=0}^{(k-1)N -1}\frac{1}{\omega_{j,p}} \mathbb{P}_{(k-1)N} \quad \forall \quad k \geq p}$

(ii) $\displaystyle{\mathbb{H}_{p}^{N}\prod_{j=s}^{kN -1+s}\frac{1}{\omega_{j,p}}\mathbb{P}_{kN+s}} = \displaystyle{\prod_{j=s}^{(k-1)N -1}\frac{1}{\omega_{j,p}} \mathbb{P}_{(k-1)N+s}}$ for $s\geq p$,$N\geq s$ and $k \geq p$\\

(iii) $\displaystyle{\varphi_{s,N}}$ is $N$-periodic point of  $\mathbb{H}_{p}$.\\

(iv) $\displaystyle{\varphi_{s,N} \in D(\mathbb{H}_{p}^{N})}$.\\

\end{Lemma}

Now, Let \\

$\phi(z) = \displaystyle{\sum_{s=p}^{M}a_{s}\mathbb{P}_{s}(z)}$ $\hfill { }  (4.8)$\\

such that \\

$\displaystyle{\mid a_{s}\prod_{j=p}^{s-1}\omega_{j,p}\mid < 1 ; s= p, p+1, ........, M}$ $\hfill { }  (4.9)$\\

and we choose the periodic point for $\mathbb{H}_{p}$ as $\psi(z)$\\

 $\psi(z) = \displaystyle{\sum_{s=p}^{M}a_{s}\varphi_{s,N}(z)}$ $\hfill { }  (4.10)$\\

then there exists an $N \geq M$ such that\\

$\mid\mid \phi - \psi \mid\mid \leq \epsilon \quad \forall \quad \epsilon > 0$. $\hfill { }  (4.11)$\\

\end{document}